**Bessem Samet and Habib Yazidi**

# Coupled fixed point theorems in partially ordered ε-chainable metric spaces

**Bessem Samet** [1][*] **and Habib Yazidi** [2]

[1] Département de Mathématiques, Ecole Supérieure des Sciences et Techniques de Tunis
5 Avenue Taha-Hussein, B.P. :56, Bab Menara-1008, Tunisie

[2] Département de Mathématiques, Ecole Supérieure des Sciences et Techniques de Tunis
5 Avenue Taha-Hussein, B.P. :56, Bab Menara-1008, Tunisie

**ABSTRACT**

**In this paper, we introduce the notion of partially ordered ε-chainable metric spaces and we derive new coupled fixed point theorems for uniformly locally contractive mappings on such spaces.**

## 1. Introduction

The Banach fixed point theorem [4] is a simple and powerful theorem with a wide range of applications, including iterative methods for solving linear, nonlinear, differential, and integral equations. This theorem has been generalized and extended by many authors in various ways; see ([1-3], [5]-[24]) and others.

Recently, Ran and Reurings [20], Bhaskar and Lakshmikantham [9], Nieto and Lopez [18], Agarwal, El-Gebeily and O'Regan [1] and Lakshmikantham and Ciric [11] presented some new results for contractions in partially ordered metric spaces (see also [3], [5], [6], [10], [12-17], [19], [21]). For a given partially ordered set $X$, Bhaskar and Lakshmikantham in [9] introduced the concept of coupled fixed point of a mapping $F: X \times X \to X$. Later in [11] Lakshmikantham and Ciric investigated some more coupled fixed point theorems in partially ordered sets. Very recently, Samet [21] extended the results of Bhaskar and Lakshmikantham [9] to mappings satisfying a generalized Meir-Keeler contractive condition.

In this paper, we introduce the notion of partially ordered ε-chainable metric spaces and we derive new coupled fixed point theorems for uniformly locally contractive mappings on such spaces. To begin, we first recall some definitions given in [9] which will be used in this work.

**Definition 1.1.** Let $(X, \leq)$ be a partially ordered set and $F: X \times X \to X$ be a given mapping. We say that $F$ has the mixed monotone property if for any $x, y \in X$, we have:

$$x_1, x_2 \in X, x_1 \leq x_2 \Rightarrow F(x_1, y) \leq F(x_2, y)$$
$$y_1, y_2 \in X, y_1 \leq y_2 \Rightarrow F(x, y_1) \geq F(x, y_2)$$

**Definition 1.2.** Let $X$ be a non-empty set and $F : X \times X \to X$ be a given mapping. We say that $(x, y) \in X \times X$ is a coupled fixed point of $F$ if: $F(x, y) = x$ and $F(y, x) = y$.

Now, we introduce the following definitions.

**Definition 1.3.** Let $(X, \leq)$ be a partially ordered set endowed with a metric $d$ and $\varepsilon > 0$. We say that $X$ is $\varepsilon$-chainable with respect to the partial order $\leq$ on $X$, if for any two points $a, b \in X$ such that $a \leq b$, there exists a finite set of points:

$$a = \alpha_0 \leq \alpha_1 \leq \cdots \leq \alpha_{n-1} \leq \alpha_n = b$$

such that $d(\alpha_{i-1}, \alpha_i) < \varepsilon$ for all $i = 1, 2, \cdots, n$.

**Definition 1.4.** Let $(X, \leq)$ be a partially ordered set endowed with a metric $d$ and $F : X \times X \to X$ be a given mapping. We say that $F$ is $(\varepsilon, \lambda)$ uniformly locally contractive if:

$$\frac{d(x,u) + d(y,v)}{2} < \varepsilon \Rightarrow d(F(x,y), F(u,v)) < \frac{\lambda}{2}[d(x,u) + d(y,v)], \forall x \geq u, \forall y \leq v,$$

where $\varepsilon > 0$ and $\lambda \in (0,1)$.

Through this paper, we will use the following notations. Let $(X, \leq)$ be a partially ordered set endowed with a metric $d$ and $F : X \times X \to X$ be a given mapping.

- We endow the product space $X \times X$ with the partial order $\leq$ defined by:

$$(x, y), (u, v) \in X \times X, (u, v) \leq (x, y) \Leftrightarrow x \geq u, y \leq v.$$

- We endow the product space $X \times X$ with the metric $\eta$ defined by:

$$\eta((x, y), (u, v)) = d(x, u) + d(y, v), \forall (x, y), (u, v) \in X \times X.$$

- For all $(x, y) \in X \times X$, we denote:

$$F^0(x, y) = x, \; F^1(x, y) = F(x, y), \; F^{m+1}(x, y) = F(F^m(x, y), F^m(y, x)) \; \forall m \in N.$$

Here, $N$ is the set of all positive integers.

## 2. Main results

The following lemma is the principal tool used to prove the main results.

**Lemma 2.1.** Let $(X, \leq)$ be a partially ordered set endowed with a metric $d$ and $F : X \times X \to X$ be a given mapping. We assume that

1. $X$ is $\varepsilon$-chainable with respect to the partial order $\leq$ on $X$,

2. $F$ has the mixed monotone property,

3. $F$ is $(\varepsilon, \lambda)$ uniformly locally contractive mapping,

4. $\exists \ (a,b),(a^*,b^*) \in X \times X$ such that $a \leq b$ and $a^* \geq b^*$.

Then,

(1) $\quad \lim_{m \to +\infty} \eta((F^m(a,a^*), F^m(a^*,a)),(F^m(b;b^*), F^m(b^*,b))) = 0.$

Moreover, we have:

$$\eta((F^m(a,a^*), F^m(a^*,a)),(F^m(b;b^*), F^m(b^*,b))) < 2n\lambda^m, \ \forall m.$$

**Proof.** Since $X$ is $\varepsilon$-chainable, there exist $\alpha_0, \alpha_1, \cdots, \alpha_n \in X$ and $\beta_0, \beta_1, \cdots, \beta_n \in X$ such that

(2) $\quad \begin{cases} a = \alpha_0 \leq \alpha_1 \leq \cdots \leq \alpha_{n-1} \leq \alpha_n = b \\ d(\alpha_{i-1}, \alpha_i) < \varepsilon, \forall i = 1,2,\cdots,n \end{cases}$

and

(3) $\quad \begin{cases} b^* = \beta_n \leq \beta_{n-1} \leq \cdots \leq \beta_1 \leq \beta_0 = a^* \\ d(\beta_{i-1}, \beta_i) < \varepsilon, \forall i = 1,2,\cdots,n. \end{cases}$

From (2)-(3) and using the mixed monotone property of $F$, we can show easily that for all $i$, we have:

(4) $\quad F^m(\alpha_i, \beta_i) \geq F^m(\alpha_{i-1}, \beta_{i-1})$ and $F^m(\beta_i, \alpha_i) \leq F^m(\beta_{i-1}, \alpha_{i-1})$ for all $m \in N$.

Now, we claim that for all $m \in N$, we have:

(5) $\quad d(F^m(\alpha_i, \beta_i), F^m(\alpha_{i-1}, \beta_{i-1})) < \lambda^m \varepsilon$ and $d(F^m(\beta_i, \alpha_i), F^m(\beta_{i-1}, \alpha_{i-1})) < \lambda^m \varepsilon.$

To prove (5), we will argue by induction. This result is trivial for $m = 0$. Let us check that (5) is true for $m = 1$. From (2)-(3), we get:

$$\frac{d(\alpha_i, \alpha_{i-1}) + d(\beta_i, \beta_{i-1})}{2} < \varepsilon, \quad \frac{d(\beta_{i-1}, \beta_i) + d(\alpha_{i-1}, \alpha_i)}{2} < \varepsilon.$$

Since $\alpha_i \geq \alpha_{i-1}$, $\beta_i \leq \beta_{i-1}$ and $F$ is $(\varepsilon, \lambda)$ uniformly locally contractive, we obtain:

$$d(F(\alpha_i, \beta_i), F(\alpha_{i-1}, \beta_{i-1})) < \lambda \varepsilon \text{ and } d(F(\beta_i, \alpha_i), F(\beta_{i-1}, \alpha_{i-1})) < \lambda \varepsilon.$$

Then, (5) is true for $m = 1$. Now, assume that (5) holds for a given $m \in N$. Let us prove that (5) holds also for $m + 1$.

Since (5) holds for $m$, we get:

$$\frac{d(F^m(\alpha_i,\beta_i),F^m(\alpha_{i-1},\beta_{i-1}))+d(F^m(\beta_i,\alpha_i),F^m(\beta_{i-1},\alpha_{i-1}))}{2}<\lambda^m\varepsilon$$

$$\frac{d(F^m(\beta_{i-1},\alpha_{i-1}),F^m(\beta_i,\alpha_i))+d(F^m(\alpha_{i-1},\beta_{i-1}),F^m(\alpha_i,\beta_i))}{2}<\lambda^m\varepsilon.$$

Then, from (4) and since $F$ is $(\varepsilon,\lambda)$ uniformly locally contractive, we obtain:

$$\begin{cases} d(F(F^m(\alpha_i,\beta_i),F^m(\beta_i,\alpha_i)),F(F^m(\alpha_{i-1},\beta_{i-1}),F^m(\beta_{i-1},\alpha_{i-1})))<\lambda^{m+1}\varepsilon \\ d(F(F^m(\beta_{i-1},\alpha_{i-1}),F^m(\alpha_{i-1},\beta_{i-1})),F(F^m(\beta_i,\alpha_i),F^m(\alpha_i,\beta_i)))<\lambda^{m+1}\varepsilon \end{cases}$$

which implies that (5) holds for $m+1$. Then, (5) holds for all $m \in N$. Now, using the triangular inequality and (5), we get:

(6)
$$d(F^m(a,a^*),F^m(b,b^*)) \leq d(F^m(\alpha_0,\beta_0),F^m(\alpha_1,\beta_1))+\cdots+d(F^m(\alpha_{n-1},\beta_{n-1}),F^m(\alpha_n,\beta_n))<n\lambda^m\varepsilon$$

Similarly, one can show that

(7) $\quad d(F^m(a^*,a),F^m(b^*,b))<n\lambda^m\varepsilon.$

Combining (6) and (7) and using that $\lambda \in (0,1)$, we obtain:

$$\eta((F^m(a,a^*),F^m(a^*a)),(F^m(b;b^*),F^m(b^*,b)))<2n\lambda^m \to 0 \text{ as } m \to +\infty.$$

This makes end to the proof.

Now, we are able to prove some theorems. We start by studying the existence of a coupled fixed point. Our first result is the following.

**Theorem 2.1.** Let $(X,\leq)$ be a partially ordered set endowed with a metric $d$ such that $(X,d)$ is complete. Let $F:X \times X \to X$ be a given mapping. We assume that

1. $X$ is $\varepsilon$-chainable with respect to the partial order $\leq$ on $X$,
2. $F$ is continuous,
3. $F$ has the mixed monotone property,
4. $F$ is $(\varepsilon,\lambda)$ uniformly locally contractive mapping,
5. $\exists\ x_0,y_0 \in X$ such that $x_0 \leq F(x_0,y_0)$ and $y_0 \geq F(y_0,x_0)$.

Then, $F$ admits a coupled fixed point.

**Proof.** Let us define the sequences $\{x_m\}$ and $\{y_m\}$ in $X$ by :

$$\begin{cases} x_{m+1} = F^{m+1}(x_0, y_0) = F(F^m(x_0, y_0), F^m(y_0, x_0)), \\ y_{m+1} = F^{m+1}(y_0, x_0) = F(F^m(y_0, x_0), F^m(x_0, y_0)). \end{cases}$$

By taking $(a,b) = (x_0, F(x_0, y_0)) = (x_0, x_1)$ and $(a^*, b^*) = (y_0, F(y_0, x_0)) = (y_0, y_1)$, we show that all the hypotheses required by Lemma 2.1 are satisfied. Hence,

(8) $\quad A_m := \eta((F^m(x_0, y_0), F^m(y_0, x_0)), (F^m(x_1, y_1), F^m(y_1, x_1))) < 2n\lambda^m \varepsilon.$

Now, we will show that $\{x_m\}$ and $\{y_m\}$ are Cauchy sequences in $X$. We have :

$$d(x_m, x_{m+1}) = d(F^m(x_0, y_0), F^{m+1}(x_0, y_0)) = d(F^m(x_0, y_0), F^m(x_1, y_1)) \leq A_m.$$

Then, from (8), it follows immediately that $\{x_m\}$ is a Cauchy sequence in $X$. Similarly, we have :

$$d(y_m, y_{m+1}) = d(F^m(y_0, x_0), F^{m+1}(y_0, x_0)) = d(F^m(y_0, x_0), F^m(y_1, x_1)) \leq A_m$$

and $\{y_m\}$ is also a Cauchy sequence in $X$.

Now, since $(X, d)$ is a complete metric space, there exists $(x, y) \in X \times X$ such that

(9) $\quad x_m \xrightarrow{d} x \quad \text{and} \quad y_m \xrightarrow{d} y \quad \text{as } m \to +\infty.$

With the continuity of $F$, (9) implies that

(10) $F(x_m, y_m) \xrightarrow{d} F(x, y)$ and $F(y_m, x_m) \xrightarrow{d} F(y, x)$ as $m \to +\infty.$

Using the triangular inequality, (9) and (10), we obtain:

$$d(x, F(x, y)) \leq d(x, x_{m+1}) + d(F(x_m, y_m), F(x, y)) \to 0 \text{ as } m \to +\infty.$$

This implies that $F(x, y) = x$. Similarly, we have :

$$d(y, F(y, x)) \leq d(y, y_{m+1}) + d(F(y_m, x_m), F(y, x)) \to 0 \text{ as } m \to +\infty.$$

Then, $F(y, x) = y$. Finally, $(x, y)$ is a coupled fixed point of $F$. This makes end to the proof.

As it is showed in [9], if we require that the underlying metric space $X$ has an additional property, the previous result is still valid for $F$ not necessarily continuous. We discuss this in the following theorem.

**Theorem 2.2.** Let $(X, \leq)$ be a partially ordered set endowed with a metric $d$ such that $(X, d)$ is complete. Let $F : X \times X \to X$ be a given mapping. We assume that

1. $X$ is $\varepsilon$-chainable with respect to the partial order $\leq$ on $X$,

2. if $\{x_m\}$ is a nondecreasing sequence in $X$ such that $x_m \xrightarrow{d} x$ as $m \to +\infty$, then $x_m \leq x$ for all $m$,

3. if $\{y_m\}$ is a nonincreasing sequence in $X$ such that $y_m \xrightarrow{d} y$ as $m \to +\infty$, then $y_m \geq y$ for all $m$,

4. $F$ has the mixed monotone property,

5. $F$ is $(\varepsilon, \lambda)$ uniformly locally contractive mapping,

6. $\exists\ x_0, y_0 \in X$ such that $x_0 \leq F(x_0, y_0)$ and $y_0 \geq F(y_0, x_0)$.

Then, $F$ admits a coupled fixed point.

**Proof.** Following the proof of Theorem 2.1, we have only to show that $(x, y)$ is a coupled fixed point of $F$. Let $p > 1$. From 1 and 2, there exists $m(p) \in N$ such that

$$(11)\ \begin{cases} \dfrac{d(x, x_{m(p)}) + d(y, y_{m(p)})}{2} = \dfrac{d(y_{m(p)}, y) + d(x_{m(p)}, x)}{2} < \dfrac{\varepsilon}{p} (< \varepsilon), \\ d(x, x_{m(p)+1}) < \dfrac{\varepsilon}{p}, d(y, y_{m(p)+1}) < \dfrac{\varepsilon}{p}. \end{cases}$$

From the hypothesis 4, it is clear that $\{x_m\}$ is a nondecreasing sequence and $\{y_m\}$ is a nonincreasing sequence. Then, from hypotheses 2 and 3, we have:

$$(12)\ x_m \leq x \text{ and } y_m \geq y \text{ for all } m.$$

Since $F$ is $(\varepsilon, \lambda)$ uniformly locally contractive, from (11)-(12), we get:

$$(13)\ d(F(x, y), F(x_{m(p)}, y_{m(p)})) < \dfrac{\lambda \varepsilon}{p} \text{ and } d(F(y_{m(p)}, x_{m(p)}), F(y, x)) < \dfrac{\lambda \varepsilon}{p}.$$

Using the triangular inequality, from (11) and (13), we obtain :

$$d(x, F(x, y)) \leq d(x, x_{m(p)+1}) + d(F(x_{m(p)}, y_{m(p)}), F(x, y)) < \dfrac{\varepsilon(\lambda+1)}{p} \to 0 \text{ as } p \to +\infty.$$

Hence, $x = F(x, y)$. By a similar argument, we can show that $y = F(y, x)$. Finally, $(x, y)$ is a coupled fixed point of $F$ and the proof is completed.

One can prove that the coupled fixed point is in fact unique, provided that the product space $X \times X$ endowed with the partial order mentioned earlier has the following property:

(H): $\forall (x, y), (x^*, y^*) \in X \times X, \exists (z_1, z_2) \in X \times X$ that is comparable to $(x, y)$ and $(x^*, y^*)$.

This is the purpose of the next theorem.

**Theorem 2.3.** Adding condition (H) to the hypotheses of Theorem 2.1, we obtain the uniqueness of the coupled fixed point of $F$.

**Proof.** Assume that $(x^*, y^*)$ is another coupled fixed point of $F$. We distinguish two cases.

<u>First case:</u> $(x, y)$ and $(x^*, y^*)$ are comparable with respect to the ordering in $X \times X$. Without restriction to the generality we can assume that $x \leq x^*$ and $y \geq y^*$. Applying Lemma 2.1, we get:

$$\lim_{m \to +\infty} \eta((F^m(x, y), F^m(y, x)), (F^m(x^*; y^*), F^m(y^*, x^*))) = 0$$

On the other hand, for all $m \in N$, we have:

$$x = F^m(x, y), \ y = F^m(y, x), \ x^* = F^m(x^*, y^*), \ y^* = F^m(y^*, x^*).$$

Then, $\eta((x, y), (x^*, y^*)) = 0$ and $(x, y) = (x^*, y^*)$.

<u>Second case:</u> $(x, y)$ and $(x^*, y^*)$ are not comparable. From (H), there exists $(z_1, z_2) \in X \times X$ that is comparable to $(x, y)$ and $(x^*, y^*)$. Without restriction to the generality, we can suppose that $x \leq z_1$, $y \geq z_2$ and $x^* \leq z_1$, $y^* \geq z_2$. Again, applying Lemma 2.1, we get:

$$(14) \quad \begin{cases} \lim_{m \to +\infty} \eta((F^m(x, y), F^m(y, x)), (F^m(z_1, z_2), F^m(z_2, z_1))) = 0, \\ \lim_{m \to +\infty} \eta((F^m(x^*, y^*), F^m(y^*, x^*)), (F^m(z_1, z_2), F^m(z_2, z_1))) = 0. \end{cases}$$

Now, using the triangular inequality and (14), we obtain:

$$\eta((x, y), (x^*, y^*)) = \eta((F^m(x, y), F^m(y, x)), (F^m(x^*, y^*), F^m(y^*, x^*)))$$

$$\leq \eta((F^m(x, y), F^m(y, x)), (F^m(z_1, z_2), F^m(z_2, z_1)))$$

$$+ \eta((F^m(z_1, z_2), F^m(z_2, z_1)), (F^m(x^*, y^*), F^m(y^*, x^*)))$$

$$\to 0 \text{ as } m \to +\infty.$$

Then, $\eta((x, y), (x^*, y^*)) = 0$ and $(x, y) = (x^*, y^*)$. This makes end to the proof.

Now, we will prove the following result.

**Theorem 2.4.** In addition to the hypotheses of Theorem 2.1, suppose that every pair of elements of $X$ has an upper or a lower bound in $X$. Then, $x = y$.

**Proof.** We distinguish two cases.

<u>First case:</u> $x = F^m(x, y)$ is comparable to $y = F^m(y, x)$. We can assume that $x \leq y$. We can write:
$$x \leq y \text{ and } y \geq y.$$
Applying Lemma 2.1, we get:

(15) $\lim_{m \to +\infty} B_m := \eta((F^m(x, y), F^m(y, x)), (F^m(y, y), F^m(y, y))) = 0.$

From (15), we get:
$$d(x, y) = d(F^m(x, y), F^m(y, x)) \leq d(F^m(x, y), F^m(y, y)) + d(F^m(y, y), F^m(y, x))$$
$$= B_m \to 0 \text{ as } m \to +\infty.$$

Then, $x = y$.

<u>Second case :</u> $x$ is not comparable to $y$. Then, there exists an upper bound or lower bound of $x$ and $y$. That is, there exists $z \in X$ comparable with $x$ and $y$. For example, we can suppose that $x \leq z$ and $y \leq z$. Again, applying Lemma 2.1, we obtain:

(16) $\lim_{m \to +\infty} C_m := \eta((F^m(x, y), F^m(y, x)), (F^m(z, z), F^m(z, z))) = 0.$

From (16), we get:
$$d(x, y) = d(F^m(x, y), F^m(y, x)) \leq d(F^m(x, y), F^m(z, z)) + d(F^m(z, z), F^m(y, x))$$
$$= C_m \to 0 \text{ as } m \to +\infty.$$

Then, $x = y$ and the proof is completed.

Alternatively, if we know that the elements $x_0$ and $y_0$ are such that $x_0 \leq y_0$, then we can also demonstrate that the components $x$ and $y$ of the coupled fixed point are indeed the same. This is the purpose of the next theorem.

**Theorem 2.5.** In addition to the hypotheses of Theorem 2.1 (resp. Theorem 2.2), suppose that $x_0, y_0 \in X$ are comparable. Then, $x = y$.

**Proof.** Without restriction to the generality, we can assume that $x_0 \leq y_0$. Applying Lemma 2.1, we get:

(17) $\lim_{m \to +\infty} D_m := \eta((F^m(x_0, y_0), F^m(y_0, x_0)), (F^m(y_0, x_0), F^m(x_0, y_0))) = 0.$

From (17) and using the triangular inequality, we get:

$$\begin{aligned}d(x,y) &\leq d(x,x_m)+d(x_m,y_m)+d(y_m,y)\\ &= d(x_m,x)+d(F^m(x_0,y_0),F^m(y_0,x_0))+d(y_m,y)\\ &\leq d(x_m,x)+D_m+d(y_m,y)\\ &\to 0 \text{ as } m\to+\infty.\end{aligned}$$

Then, $d(x,y)=0$ and $x=y$. This makes end to the proof.

**Acknowledgments.** The authors thank the referees for their valuable comments and suggestions.